\documentclass[11pt]{article}
\usepackage{amsmath,amssymb,cite,graphicx}

\numberwithin{equation}{section}

\newcommand{\R}{\mathbb{R}}

\newcommand{\<}{\langle}
\renewcommand{\>}{\rangle}
\newcommand{\supp}{\operatorname{supp}}
\newcommand{\sgn}{\operatorname{sgn}}
\newcommand{\E}{\operatorname{\mathbf{E}}}
\newcommand{\Prob}{\operatorname{\mathbf{P}}}

\renewcommand{\P}{\Prob}
\newcommand{\UOT}{U_{\Omega T}}
\newcommand{\UtU}{\UOT^*\UOT}
\newcommand{\ukuk}{u^k \otimes u^k}
\newcommand{\uk}{u^k}
\newcommand{\cF}{{\cal F}}
\newcommand{\cB}{\mathcal{B}}

\newenvironment{proof}{\noindent {\bf Proof} }{\endprf\par}
\def \endprf{\hfill {\vrule height6pt width6pt depth0pt}\medskip} 
\newtheorem{theorem}{Theorem}[section]

\newtheorem{lemma}{Lemma}[section]

\begin{document}

\title{Sparsity and Incoherence in Compressive Sampling}

\author{Emmanuel Cand\`es$^{\dagger}$ and Justin Romberg$^{\sharp}$ \\
\vspace{-.2cm}\\
$\dagger$ Applied and Computational Mathematics, Caltech, Pasadena, CA 91125 \\
\vspace{-.3cm}\\
$\sharp$ Electrical and Computer Engineering, Georgia Tech, Atlanta, GA  90332}

\date{November 2006}

\maketitle

\begin{abstract}
  We consider the problem of reconstructing a sparse signal
  $x^0\in\R^n$ from a limited number of linear measurements.  Given
  $m$ randomly selected samples of $U x^0$, where $U$ is an
  orthonormal matrix, we show that $\ell_1$ minimization recovers
  $x^0$ exactly when the number of measurements exceeds
\[
m\geq \mathrm{Const}\cdot\mu^2(U)\cdot S\cdot\log n,
\]
where $S$ is the number of nonzero components in $x^0$, and $\mu$ is
the largest entry in $U$ properly normalized: $\mu(U) = \sqrt{n} \cdot
\max_{k,j}~ |U_{k,j}|$. The smaller $\mu$, the fewer samples needed. 

The result holds for ``most'' sparse signals $x^0$ supported on a
fixed (but arbitrary) set $T$.  Given $T$, if the sign of $x^0$ for
each nonzero entry on $T$ and the observed values of $Ux^0$ are drawn
at random, the signal is recovered with overwhelming
probability. Moreover, there is a sense in which this is nearly
optimal since any method succeeding with the same probability would
require just about this many samples.
\end{abstract}

%\vfill
{\small {\bf Acknowledgments.} E.~C. is partially supported by
  National Science Foundation grants ITR ACI-0204932 and CCF–515362,
  by the 2006 Waterman Award (NSF), and by a grant from
  DARPA. J.~R.~is supported by National Science Foundation grants
  CCF–515362 and ITR ACI-0204932. E.~C. would like to thank Joel Tropp
  and Houman Owhadi for fruitful conversations related to this
  project.

}
%\pagebreak

%---------------------------------------------------------------------------------
\section{Introduction}

\subsection{Sparse recovery from partial measurements}

This paper addresses the problem of signal acquisition in a broad setting.  
We are interested in ``sampling'' a vector
$x^0\in\R^n$.  Instead of observing $x^0$ directly, we sample a small
number $m$ of {\em transform coefficients} of $x^0$.  For an
orthogonal matrix\footnote{On a first reading, our choice of
  normalization of $U$ may seem a bit strange.  The advantages of
  taking the row vectors of $U$ to have Euclidean norm $\sqrt{n}$ are
  that 1) the notation in the sequel will be cleaner, and 2) it will
  be easier to see how this result generalizes the special case of
  incomplete sampling in the Fourier domain presented in
  \cite{candes06ro}.} $U$ with
\begin{equation}
\label{eq:UtU}
U^*U = nI,
\end{equation}
these transform coefficients are given by $y^0 =U x^0$.  Of course, if
all $n$ of the coefficients $y^0$ are observed, recovering $x^0$ is
trivial: we simply apply $\frac{1}{n} \, U^*$ to the vector of
observations $y^0$.  Instead, we are concerned with the highly
underdetermined case in which only a small fraction of the components
of $y^0$ are actually sampled or observed.  Given a subset
$\Omega\subset\{1,\ldots,n\}$ of size $|\Omega|=m$, the challenge is
to infer the ``long'' $n$-dimensional vector $x^0$ from the ``short''
$m$-dimensional vector of observations $y=U_\Omega x^0$, where
$U_\Omega$ is the $m\times n$ matrix consisting of the rows of $U$
indexed by $\Omega$. In plain English, we wish to solve a system of
linear equations in which there are fewer equations than unknowns.

A special instance of this problem was investigated in a recent paper
\cite{candes06ro}, where $U$ is taken as the usual discrete Fourier
transform.  The main result of this work is that if $x^0$ is {\em
  $S$-sparse} (at most $S$ of the $n$ components of $x^0$ are
nonzero), then it can be recovered perfectly from on the order of
$S\log n$ Fourier-domain samples.  The recovery algorithm is concrete
and tractable: given the discrete Fourier coefficients
\begin{equation}
  \label{eq:dft}
  y_k = \sum_{t=1}^{n} 
x^0(t) e^{-i2\pi (t-1)k/n},~~k \in \Omega,   
\end{equation}
or $y = F_\Omega x^0$ for short, we solve the convex optimization
program
\[
\min_x~\|x\|_{\ell_1}\quad\text{subject to}\quad F_\Omega x = y.
\]
For a fixed $x^0$, the recovery is exact for the overwhelming majority of
sample sets $\Omega$ of size
\begin{equation}
\label{eq:fourierrec}
|\Omega| \geq C\cdot S\cdot\log n,
\end{equation}
where $C$ is a known (small) constant.  

%%%%

Since \cite{candes06ro}, a theory of ``compressed sensing'' has
developed around several papers
\cite{candes04ne,candes05de,donoho04co} demonstrating the
effectiveness of $\ell_1$ minimization for recovering sparse signals
from a limited number of measurements.  To date, most of this effort
has been focused on systems which take completely unstructured,
noise-like measurements, i.e.\ the observation vector $y$ is created
from a series of inner products against random test vectors
$\{\phi_k\}$:
\begin{equation}
\label{eq:randomy}
y_k = \<\phi_k,x^0\>,\quad k=1,\ldots,m.
\end{equation}
The collection $\{\phi_k\}$ is sometimes referred to as a {\em
  measurement ensemble}; we can write \eqref{eq:randomy} compactly as
$y=\Phi x^0$, where the rows of $\Phi$ are the $\phi_k$.  Published
results take $\phi_k$ to be a realization of Gaussian white noise, or
a sequence of Bernoulli random variables taking values $\pm 1$ with
equal probability.  This work has shown that taking random
measurements is in some sense an optimal strategy for acquiring sparse
signals; it requires a near-minimal number of measurements
\cite{candes04ne,candes05de,donoho04co,donoho04fo,baraniuk06ra} ---
$m$ measurements can recover signals with sparsity $S\lesssim
m/\log(n/m)$, and all of the constants appearing in the analysis are
small \cite{donoho05ne}.
Similar bounds have also appeared using greedy \cite{tropp05si} and complexity-based \cite{haupt06si} recovery algorithms in place of $\ell_1$ minimization.

Although theoretically powerful, the practical relevance of results
for completely random measurements is limited in two ways.
The first is that we are not always at liberty to choose the types of
measurements we use to acquire a signal.  For example, in magnetic
resonance imaging (MRI), subtle physical properties of nuclei are
exploited to collect samples in the Fourier domain of a two- or
three-dimensional object of interest.  While we have control over
which Fourier coefficients are sampled, the measurements are
inherently frequency based.  A similar statement can be made about
tomographic imaging; the machinery in place measures Radon slices, and
these are what we must use to reconstruct an image.

The second drawback to completely unstructured measurement systems is
computational.  Random (i.e.\ unstructured) measurement ensembles are
unwieldy numerically; for large values of $m$ and $n$, simply storing
or applying $\Phi$ (tasks which are necessary to solve the $\ell_1$
minimization program) are nearly impossible.
If, for example, we want to reconstruct a megapixel image
($n=1,000,000$) from $m=25,000$ measurements (see the numerical
experiment in Section~\ref{sec:applications}), we would need more than
3 gigabytes of memory just to store the measurement matrix, and on the
order of gigaflops to apply it.
The goal from this point of view, then, is to have similar recovery
bounds for measurement matrices $\Phi$ which can be applied quickly
(in $O(n)$ or $O(n\log n)$ time) and implicitly (allowing us to use a
``matrix free'' recovery algorithm).

Our main theorem, stated precisely in Section~\ref{sec:result} and
proven in Section~\ref{sec:proof}, states that bounds analogous to
\eqref{eq:fourierrec} hold for sampling with general orthogonal
systems.  We will show that for a fixed signal support $T$ of size
$|T|=S$, the program
\begin{equation}
\label{eq:genl1}
\min_x~\|x\|_{\ell_1} \quad\text{subject to}\quad U_\Omega x = U_\Omega x^0
\end{equation}
recovers the overwhelming majority of $x^0$ supported on $T$ and observation
subsets $\Omega$ of size
\begin{equation}
\label{eq:genrec}
|\Omega| \geq C\cdot\mu^2(U)\cdot S\cdot\log n,
\end{equation}
where $\mu(U)$ is simply the largest magnitude among the
entries in $U$:
\begin{equation}
\label{eq:mu}
\mu(U) = \max_{k,j}~|U_{k,j}|.
\end{equation}

It is important to understand the relevance of the parameter $\mu(U)$
in \eqref{eq:genrec}.  $\mu(U)$ can be interpreted as a rough measure
of how concentrated the rows of $U$ are.  Since each row (or column)
of $U$ necessarily has an $\ell_2$-norm equal to $\sqrt{n}$, $\mu$
will take a value between $1$ and $\sqrt{n}$.  When the rows of $U$
are perfectly flat --- $|U_{k,j}| = 1$ for each $k,j$, as in the case
when $U$ is the discrete Fourier transform, we will have $\mu(U)=1$,
and \eqref{eq:genrec} is essentially as good as \eqref{eq:fourierrec}.
If a row of $U$ is maximally concentrated --- all the row entries but
one vanish --- then $\mu^2(U)=n$, and \eqref{eq:genrec} offers us no
guarantees for recovery from a limited number of samples. This result
is very intuitive. Suppose indeed that $U_{k_0,j_0} = \sqrt{n}$ and
$x^0$ is 1-sparse with a nonzero entry in the $j_0$th location. To
reconstruct $x^0$, we need to observe the $k_0$th entry of $Ux^0$ as
otherwise, the data vector $y$ will vanish. In other words, to
reconstruct $x^0$ with probability greater than $1 - 1/n$, we will
need to see all the components of $Ux^0$, which is just about the
content of \eqref{eq:genrec}. This shows informally that
\eqref{eq:genrec} is fairly tight on both ends of the range of the
parameter $\mu$.

For a particular application, $U$ can be decomposed as a product 
of a {\em sparsity basis} $\Psi$, and an orthogonal measurement system
$\Phi$. Suppose for instance that we wish to recover a signal $f \in
\R^n$ from $m$ measurements of the form $y = \Phi f$. The signal may
not be sparse in the time domain but its expansion in the basis $\Psi$
may be
\[
f(t) = \sum_{j = 1}^n x^0_j \psi_j(t),  \quad f = \Psi x  
\]
(the columns of $\Psi$ are the discrete waveforms $\psi_j$).  Our
program searches for the coefficient sequence in the $\Psi$-domain
with minimum $\ell_1$ norm that explains the samples in the
measurement domain $\Phi$. In short, it solves \eqref{eq:genrec} with
\[
U=\Phi\Psi,\quad \Psi^*\Psi = I,\quad \Phi^*\Phi = nI.
\]

The result \eqref{eq:genrec} then tells us how the relationship
between the sensing modality ($\Phi$) and signal model ($\Psi$)
affects the number of measurements required to reconstruct a sparse
signal.  The parameter $\mu$ can be rewritten as
\[
\mu(\Phi\Psi) = \max_{k,j}|\<\phi_k,\psi_j\>|,
\]
and serves as a rough characterization of the degree of similarity
between the sparsity and measurement systems.  For $\mu$ to be close
to its minimum value of $1$, each of the measurement vectors (rows of
$\Phi$) must be ``spread out'' in the $\Psi$ domain.  To emphasize
this relationship, $\mu(U)$ is often referred to as the {\em mutual
  coherence} \cite{donoho01un,donoho03op}.
The bound \eqref{eq:genrec} tells us that an $S$-sparse signal can be
reconstructed from $\sim S\log n$ samples in any domain in which the
test vectors are ``flat'', i.e.~the coherence parameter is $O(1)$.

%----------------------------------------------------------------------------------
\subsection{Main result}
\label{sec:result}

The ability of the $\ell_1$-minimization program \eqref{eq:genl1} to
recover a given signal $x^0$ depends only on 1) the support set $T$ of
$x^0$, and 2) the sign sequence $z_0$ of $x^0$ on $T$.\footnote{In
  other words, the recoverability of $x^0$ is determined by the {\em
    facet} of the $\ell_1$ ball of radius $\|x^0\|_{\ell_1}$ on which
  $x^0$ resides.}  For a fixed support $T$, our main theorem shows
that perfect recovery is achieved for the overwhelming majority of the
combinations of sign sequences on $T$, and sample locations (in the
$U$ domain) of size $m$ obeying \eqref{eq:genrec}.

The language ``overwhelming majority'' is made precise by introducing
a probability model on the set $\Omega$ and the sign sequence $z$.
The model is simple: select $\Omega$ uniformly at random from the set
of all subsets of the given size $m$; choose each $z(t),~t\in T$ to be
$\pm 1$ with probability $1/2$.  Our main result is:
\begin{theorem}
\label{th:main}
Let $U$ be an $n\times n$ orthogonal matrix ($U^*U=nI$) with
$|U_{k,j}|\leq\mu(U)$.  Fix a subset $T$ of the signal domain.  Choose
a subset $\Omega$ of the measurement domain of size $|\Omega|=m$, and
a sign sequence $z$ on $T$ uniformly at random. Suppose that
\begin{equation}
\label{eq:main}
m\geq C_0 \cdot |T| \cdot\mu^2(U) \cdot \log (n/\delta)
\end{equation}
and also $m \ge C'_0 \cdot \log^2(n/\delta)$ for some fixed numerical
constants $C_0$ and $C'_0$. Then with probability exceeding
$1-\delta$, every signal $x^0$ supported on $T$ with signs matching
$z$ can be recovered from $y=U_\Omega x^0$ by solving
\eqref{eq:genl1}.
\end{theorem}

The hinge of Theorem~\ref{th:main} is a new {\em weak uncertainty
  principle} for general orthobases.  Given $T$ and $\Omega$ as above,
it is impossible to find a signal which is concentrated on $T$ and on
$\Omega$ in the $U$ domain. In the example above where $U = \Phi
\Psi$, this says that one cannot be concentrated on small sets in the
$\Psi$ and $\Phi$ domains simultaneously. As noted in previous
publications \cite{candes06qu,candes06ro}, this is a statement about
the eigenvalues of minors of the matrix $U$.  Let $U_T$ be the
$n\times |T|$ matrix corresponding to the columns of $U$ indexed by
$T$, and let $U_{\Omega T}$ be the $m\times |T|$ matrix corresponding
to the rows of $U_T$ indexed by $\Omega$.  In Section~\ref{sec:proof}, 
we will prove the following:
\begin{theorem}
\label{teo:largedev}
Let $U,T$,and $\Omega$ be as in Theorem~\ref{th:main}. Suppose
that the number of measurements $m$ obeys
  \begin{equation}
    \label{eq:good}
    m \ge |T| \cdot \mu^2(U) \cdot \max(C_1\log |T|,~C_2\log(3/\delta)), 
  \end{equation}
for some positive constants $C_1,C_2$. Then 
  \begin{equation}
    \label{eq:largedev}
    \Prob \left(\Vert \frac{1}{m} \UtU - I\Vert \ge 1/2 \right) \le \delta, 
  \end{equation}
  where $\|\cdot\|$ is the standard operator $\ell_2$ norm---here, the
  largest eigenvalue (in absolute value).
\end{theorem}
For small values of $\delta$, the eigenvalues of $U^*_{\Omega
  T}U_{\Omega T}$ are all close to $m$ with high probability.  To see
that this is an uncertainty principle, let $x \in \R^n$ be a sequence
supported on $T$, and suppose that $\|m^{-1}\UOT^*\UOT-I\|\leq 1/2$.
It follows that
\begin{equation}
\label{eq:wup}
\frac{m}{2}\|x\|_{\ell_2}^2 \leq \|U_{\Omega} x\|^2_{\ell_2} 
\leq \frac{3m}{2}\|x\|^2_{\ell_2}, 
\end{equation}
which asserts that only a small portion of the energy of $x$ will be
concentrated on the set $\Omega$ in the $U$-domain (the total energy
obeys $\|Ux\|_{\ell_2}^2 = n\|x\|_{\ell_2}^2$). Moreover, this
portion is essentially proportional to the size of $\Omega$.

%----------------------------------------------------------------------------------
\subsection{Contributions and relationship to prior work}

The relationship of the mutual incoherence parameter $\mu$ to the
performance of $\ell_1$ minimization programs with equality
constraints first appeared in the context of {\em Basis Pursuit} for
sparse approximation, see \cite{donoho01un} and also
\cite{elad02ge,donoho03op,gribonval03sp}.

As mentioned in the previous section, \cite{candes06ro} demonstrated
the effectiveness of $\ell_1$ recovery from Fourier-domain samples in
slightly more general situations than in Theorem~\ref{th:main}
(randomization of the signs on $T$ is not required).  Obviously, the
results presented in this paper considerably extend this Fourier
sampling theorem.  

We also note that since \cite{candes06ro}, several
papers have appeared on using $\ell_1$ minimization to recover sparse
signals from a limited number of measurements
\cite{candes04ne,donoho04co,candes06st}.  In particular
\cite{candes04ne} and \cite{rudelson06sp} provide bounds for
reconstruction from a random subset of measurements selected from an
orthogonal basis; these papers ask that {\em all} sparse signals to be
simultaneously recoverable from the same set of samples (which is
stronger than our goal here), and their bounds have $\log$ factors of
$(\log n)^6$ and $(\log n)^5$ respectively.  These results are based on
{\em uniform} uncertainty principles, which require \eqref{eq:largedev} 
to hold for all sets $T$ of a certain size simultaneously once $\Omega$ 
is chosen.  Whether or not this $\log$ power can be reduced in this context 
remains an open question.

A contribution of
this paper is to show that if one is only interested in the recovery
of nearly all signals on a fixed set $T$, 
these extra log factors can indeed be removed. 
We show that to guarantee exact recovery, we only require $U_{\Omega T}$ to be 
well behaved for this fixed $T$ as opposed to all $T$'s of the same size,
which is a significantly weaker requirement. By examining the singular
values of $U_{\Omega T}$, one can check whether or not
$\eqref{eq:wup}$ holds.

Our method of proof, as the reader will see in
Section~\ref{sec:proof}, relies on a variation of the powerful results
presented in \cite{rudelson99ra} about the expected spectral norm of
certain random matrices.  We also introduce a novel large-deviation
inequality, similar in spirit to those reviewed in
\cite{ledoux01co,ledoux91pr} but carefully tailored for our purposes,
to turn this statement about expectation into one about high
probability.

%%%
Finally, we would like to contrast this work with \cite{tropp06ra},
which also draws on the results from \cite{rudelson99ra}.  First,
there is a difference in how the problem is framed.  In
\cite{tropp06ra}, the $m\times n$ measurement system is fixed, and
bounds for perfect recovery are derived when the support and sign
sequence are chosen at random, i.e.\ a fixed measurement system works
for most signal supports of a certain size.  In this paper, we fix an
arbitrary signal support, and show that we will be able to recover
from most sets of measurements of a certain size in a fixed domain.
Second, although slightly more general class of measurement systems is
considered in \cite{tropp06ra}, the final bounds for sparse recovery
in the context of \eqref{eq:genl1} do not fundamentally improve on the
uniform bounds cited above; \cite{tropp06ra} draws weaker conclusions
since the results are not shown to be universal in the sense that all
sparse signals are recovered as in \cite{candes04ne} and
\cite{rudelson06sp}.
%Our ``concentration of measure'' principles allow us to derive much
%tighter bounds --- in \cite{tropp06ra}, the exponent for the $(\log
%n)^\beta$ term depends on the probability of success, which is
%$O(n^{-\beta})$.

%---------------------------------------------------------------------------------
\section{Applications}
\label{sec:applications}

In the 1990s, image compression algorithms were revolutionized by the
introduction of the wavelet transform.  The reasons for this can be
summarized with two major points: the wavelet transform is a much
sparser representation for photograph-like images than traditional
Fourier-based representations, and it can be applied and inverted in
$O(n)$ computations.

To exploit this wavelet-domain sparsity in acquisition, we must have a
measurement system which is incoherent with the wavelet representation
(so that $\mu$ in \eqref{eq:genrec} is small) and that can be applied
quickly and implicitly (so that large-scale recovery is
computationally feasible).  In this section, we present numerical
experiments for two such measurement strategies.

\subsection{Fourier sampling of sparse wavelet subbands}

Our first measurement strategy takes advantage of the fact that at
fine scales, wavelets are very much spread out in frequency.  We will
illustrate this in 1D; the ideas are readily applied to 2D image.

\begin{figure}
\centerline{
\includegraphics[height=2in]{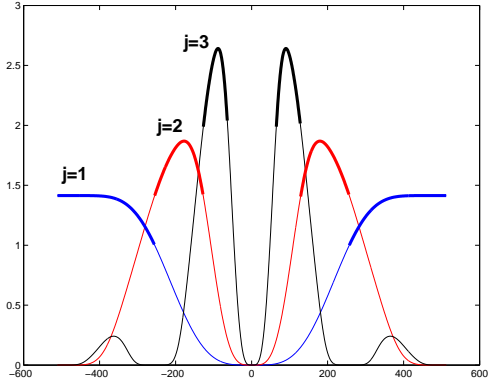}
}
\caption{\small\sl Wavelets in the frequency domain.  The curves shown above are the magnitude of the discrete Fourier transform \eqref{eq:ftwavelet} of Daubechies-8 wavelets for $n=1024$ and $j=1,2,3$.  The magnitude of $\hat{\psi}_{j,\cdot}$ over the subband \eqref{eq:sbindex} is shown in bold.}
\label{fig:waveletspectrum}
\end{figure}

Labeling the scales of the wavelet transform by $j=1,2,\ldots,J$,
where $j=1$ is the finest scale and $j=J$ the coarsest, the
wavelets\footnote{Wavelets are naturally parameterized by a scale $j$
  and a shift $k$ with $k=1,2,\ldots,n2^{-j}$ --- see
  \cite{mallat99wa}.  The wavelets at a set scale are just circular
  shifts of one another: $\psi_{j,k}(t) = \psi_{j,1}(t-2^jk)$, where
  the substraction is modulo $n$.}  $\psi_{j,k}$ at scale $j$ are
almost flat in the Fourier domain over a band of size $n_j = n2^{-j}$.
The magnitude of the Fourier transform
\begin{equation}
\label{eq:ftwavelet}
\hat{\psi}_{j,k}(\omega) = \sum_{t=1}^{n} \psi(t) e^{-i2\pi(t-1)\omega/n},
\quad \omega = -n/2+1,\ldots,n/2,
\end{equation}
is the same for each wavelet at scale $j$, since
\begin{equation}
\label{eq:wavefour}
\hat{\psi}_{j,k}(\omega) = e^{-i2\pi (k-1)\omega/n_j}\hat{\psi}_{j,1}(\omega).
\end{equation}
These spectrum magnitudes are shown for the Daubechies-8 wavelet in
Figure~\ref{fig:waveletspectrum}.  We see that over frequencies in the
$j$th subband
\begin{equation}
\label{eq:sbindex}
\omega\in\cB_j := \{n_j/2+1,\ldots,n_j\}\cup
\{-n_j+1,\ldots, -n_j/2\},
\end{equation}  
we have
\[
\frac{\max_{\omega\in\cB_j} |\hat{\psi}_{j,k}(\omega)|}{\min_{\omega\in\cB_j} |\hat{\psi}_{j,k}(\omega)|} < \mathrm{Const}\approx \sqrt{2}.
\]

%%%%%
%%%%%

Suppose now that a signal $x^0$ is a superposition of $S$ wavelets at
scale $j$, that is, we can write
\[
x^0 = \Psi_j w^0
\]
where $w^0\in\R^{n_j}$ is $S$-sparse, and $\Psi_j$ is the $n\times
n_j$ matrix whose columns are the $\psi_{j,k}(t)$ for
$k=1,\ldots,n_j$.
We will measure $x^0$ by selecting Fourier coefficients from the band
$\cB_j$ at random.  To see how this scheme fits into the domain of the
results in the introduction, let $\omega$ index the subband $\cB_j$,
let $F_j$ be the $n_j\times n$ matrix whose rows are the Fourier
vectors for frequencies in $\cB_j$, let $D_j$ be a diagonal matrix
with
\[
(D_j)_{\omega,\omega} = \hat{\psi}_{j,1}(\omega), \quad \omega \in \cB_j, 
\]
and consider the $n_j\times n_j$ system
\[
U = D^{-1}_j F_j\Psi_j.
\]
The columns of $F_j\Psi_j$ are just the Fourier transforms of the
wavelets given in \eqref{eq:wavefour},
\[
(F_j\Psi_j)_{\omega,k} = e^{-i2\pi (k-1)
  \omega/n_j}\hat{\psi}_{j,1}(\omega) \quad \Rightarrow \quad
U_{\omega,k} = e^{-i2\pi (k-1) \omega/n_j}, 
% = e^{-i\pi(k-1)}e^{i2\pi
%  (k-1)\omega'/n_j}\hat{\psi}_{j,1}(\omega),
\]
and so $U$ is just a $n_j\times n_j$ Fourier system. In fact, one can
easily check that $U^* U = U^* U = n_j \, I$.

We choose a set of Fouier coefficients $\Omega$ of size $m$ in the
band $\cB_j$, and measure
\[
y = F_\Omega x_0 = F_\Omega\Psi_j w^0,
\]
which can easily be turned into a set of samples in the $U$ domain $y'
= U_\Omega w^0$ just by re-weighting $y$.  Since the mutual
incoherence of $D^{-1}F_j\Psi_j$ is $\mu=1$, we can recover $w^0$ from
$\sim S\log n$ samples.

Table~\ref{tab:recresults} summarizes
the results of the following experiment: Fix the scale $j$, sparsity
$S$, and a number of measurements $m$.  Perform a trial for $(S,j,m)$
by first generating a signal support $T$ of size $S$, a sign sequence on that
support, and a measurement set $\Omega_j$ of size $m$ uniformly at
random, and then measuring $y=F_{\Omega_j}\Psi_j x^0$ ($x^0$ is just the sign sequence on $T$ and zero elsewhere), solving \eqref{eq:genl1}, and declaring
success if the solution matches $x^0$.  
A thousand trials were
performed for each $(S,j,m)$.  The value $M(S,j)$ recorded in the
table is the smallest value of $m$ such that the recovery was
successful in all $1000$ trials.
As with the partial Fourier ensemble (see the numerical results in
\cite{candes06ro}), we can recover from $m\approx 2S$ to $3S$
measurements.

%%%%%
\begin{table}
  \caption{\small\sl 
    Number of measurements required to reconstruct a sparse subband.  
    Here, $n=1024$, $S$ is the sparsity of the subband, and $M(S,j)$ 
    is the smallest number of measurements so that the $S$-sparse subband 
    at wavelet level $j$ was recovered perfectly in $1000/1000$ trials.}
\vspace{2mm}
\centerline{
\begin{tabular}{|c|c||c|c||c|c|}\hline
\multicolumn{2}{|c||}{$j=1$} & \multicolumn{2}{c||}{$j=2$} & \multicolumn{2}{c|}{$j=3$} \\\hline
$S$ & $M(S,j)$ & $S$ & $M(S,j)$ & $S$ & $M(S,j)$ \\\hline
50 & 100 & 25 & 56 & 15 & 35 \\\hline
25 & 68 & 15 & 40 & 8  & 24 \\\hline
15 & 49 & 8 & 27 & - & - \\\hline
\end{tabular}
}
\label{tab:recresults}
\end{table}
%%%%%

% imaging: subband-by-subband reconstruction
To use the above results in an imaging system, we would first separate
the signal/image into wavelet subband, measure Fourier coefficients in
each subband as above, then reconstruct each subband independently.
In other words, if $P_{W_j}$ is the projection operator onto the space
spanned by the columns of $\Psi_j$, we measure
\[
y^j = F_{\Omega_j} P_{W_j} x^0
\]
for $j=1,\ldots,J$, then set $w^j$ to be the solution to
\[
\min~\|w\|_{\ell_1}\quad\text{subject to}\quad F_{\Omega_j}\Psi_jw = y^j.
\]
If all of the wavelet subbands of the object we are imaging are
appropriately sparse, we will be able to recover the image perfectly.

Finally, we would like to note that this projection onto $W_j$ in the
measurement process can be avoided by constructing the wavelet and
sampling systems a little more carefully.  In \cite{donoho01un}, a
``bi-sinusoidal'' measurement system is introduced which complements
the orthonormal Meyer wavelet transform.  These bi-sinusoids are an
alternative orthobasis to the $W_j$ spanned by Meyer wavelets at a
given scale (with perfect mutual incoherence), so sampling in the
bi-sinusoidal basis isolates a given wavelet subband automatically.

In the next section, we examine an orthogonal measurement system which
allows us to forgo this subband separation all together.

\subsection{Noiselet measurements} 

In \cite{coifman01no}, a complex ``noiselet'' system is constructed
that is perfectly incoherent with the Haar wavelet representation.  If
$\Psi$ is an orthonormal system of Haar wavelets, and $\Phi$ is the
orthogonal noiselet system (renormalized so that $\Phi^*\Phi = nI$),
then $U=\Phi\Psi$ has entries of constant magnitude:
\[
|U_{k,j}| = 1,~~\forall k,j\quad\text{which
  implies}\quad \mu(U) = 1.
\]
Just as the canonical basis is maximally incoherent with the Fourier
basis, so is the noiselet system with Haar wavelets.  Thus if an
$n$-pixel image is $S$-sparse in the Haar wavelet domain, it can be
recovered (with high probability) from $\sim S\log n$ randomly
selected noiselet coefficients. 

In addition to perfect incoherence with the Haar transform, noiselets
have two additional properties that make them ideal for coded image
acquisition:
\begin{enumerate}
\item The noiselet matrix $\Phi$ can be decomposed as a multiscale
  filterbank.  As a result, it can be applied $O(n\log n)$ time.
\item The real and imaginary parts of each noiselet function are
  binary valued.  A noiselet measurement of an image is just an inner
  product with a sign pattern, which make their implementation in an
  actual acquisition system easier.  (It would be straightforward to
  use them in the imaging architecture proposed in \cite{takhar06ne},
  for example.)
\end{enumerate}

\begin{figure}
\centerline{
\begin{tabular}{ccc}
\includegraphics[height=2in]{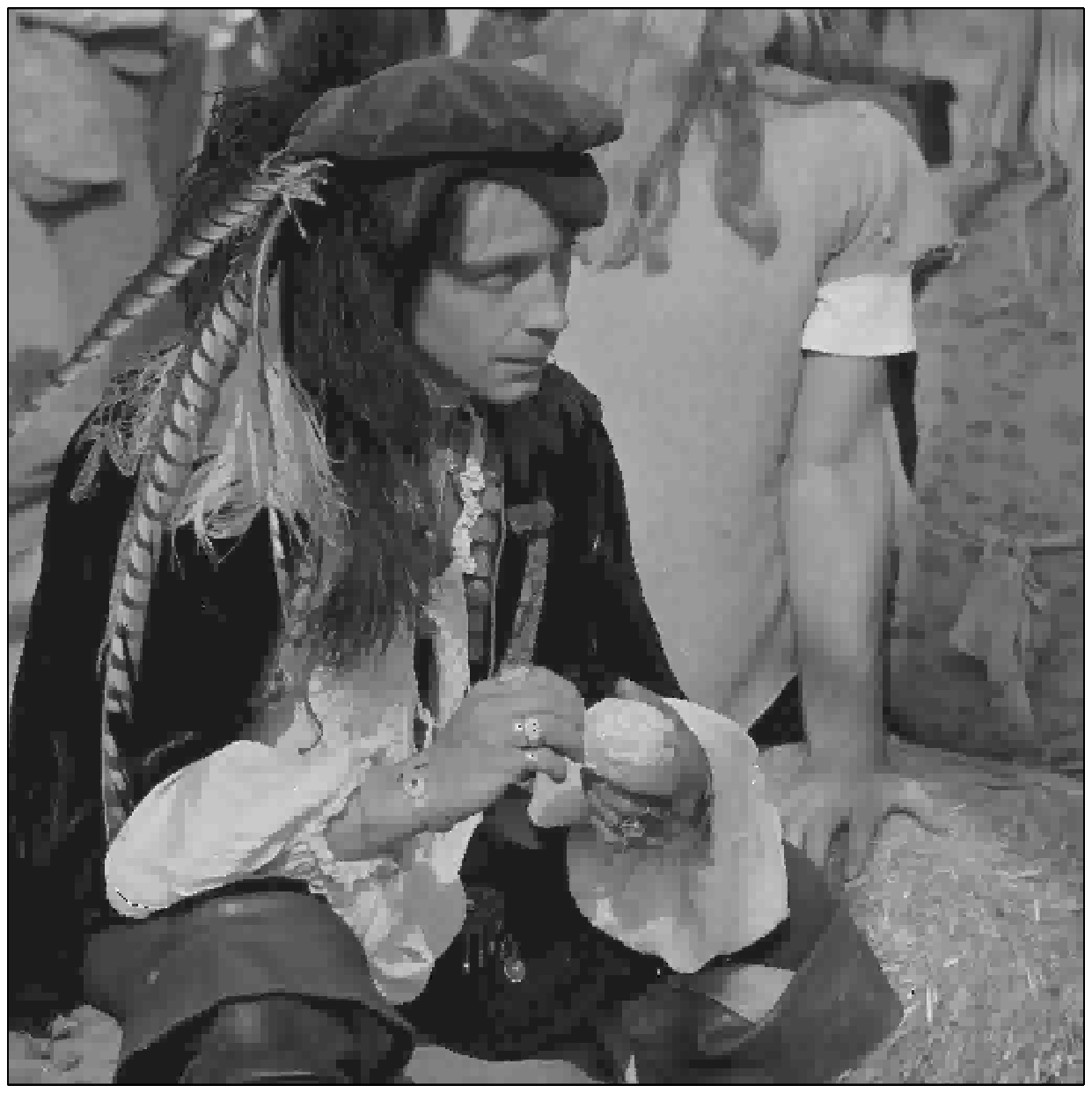} &
\includegraphics[height=2in]{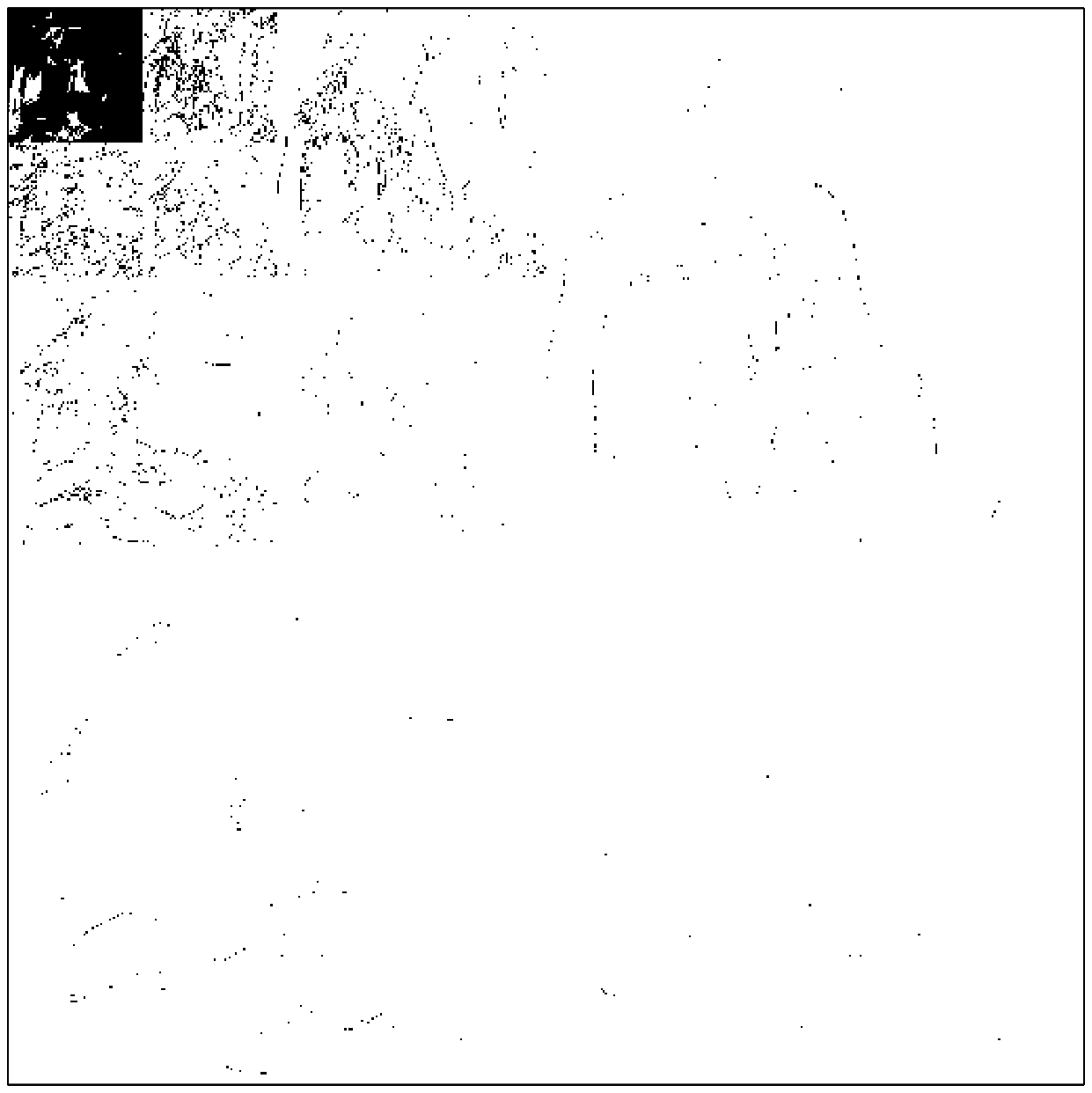} &
\includegraphics[height=2in]{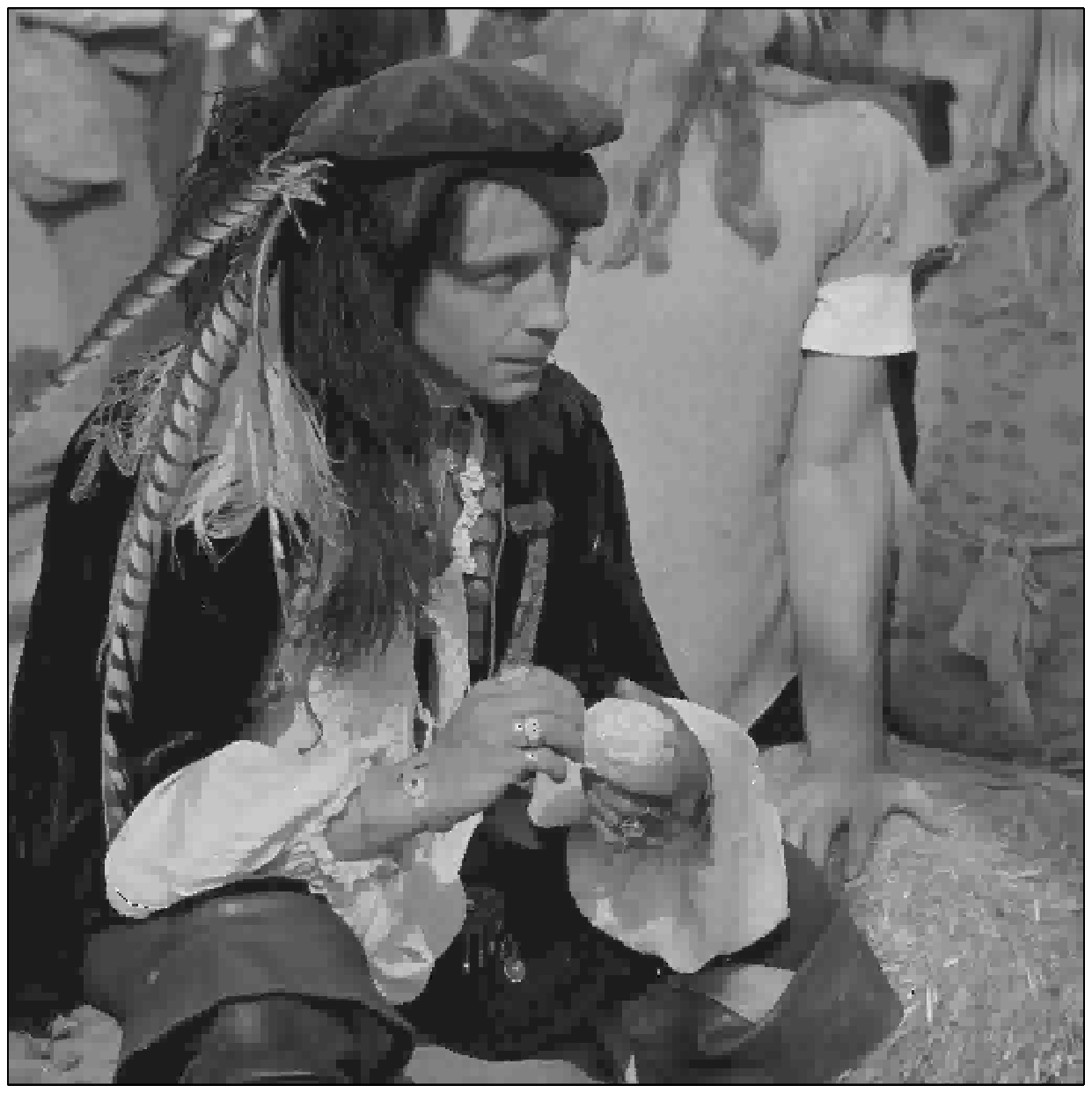} \\
(a) & (b) & (c)
\end{tabular}
}
\caption{Sparse image recovery from noiselet measurements.  (a)
  Synthetic $n=1024^2$-pixel image with $S=25,000$ non-zero Haar
  wavelet coefficients.  (b) Locations (in the wavelet quadtree) of
  significant wavelet coefficients.  (c) Image recovered from
  $m=70,000$ complex noiselet measurements.  The recovery matches (a)
  exactly.}
\label{fig:man}
\end{figure}

A large-scale numerical example is shown in Figure~\ref{fig:man}.  The
$n=1024^2$ pixel synthetic image in panel (a) is an exact
superposition of $S=25,000$ Haar wavelets\footnote{The image was
  created in the obvious way: the well-known test image was
  transformed into the Haar domain, all but the $25,000$ largest Haar
  coefficients were set to zero, and the result inverse transformed
  back into the spatial domain.}.  The observation vector $y$ was
created from $m=70,000$ randomly chosen noiselet coefficients (each
noiselet coefficient has a real and imaginary part, so there are
really $140,000$ real numbers recorded). From $y$, we are able to
recover the image exactly by solving \eqref{eq:genl1}.

This result is a nice demonstration of the compressed sensing
paradigm.  A traditional acquisition process would measure all $n\sim
10^6$ pixels, transform into the wavelet domain, and record the $S$
that are important.  Many measurements are made, but comparably very
few numbers are recorded.  Here we take only a fraction of the number
of measurements, and are able to find the $S$ active wavelets
coefficients without any prior knowledge of their locations.

The measurement process can be adjusted slightly in a practical
setting.  We know that almost all of the coarse-scale wavelet
coefficients will be important (see Figure~\ref{fig:man}(b)), so we
can potentially reduce the number of measurements needed for perfect
recovery by measuring these directly.  In fact, if we measure the
$128\times 128$ block of coarse wavelet coefficients for the image in
Figure~\ref{fig:man} directly (equivalent to measuring averages over
$8\times 8$ blocks of pixels, $16,384$ measurement total), we are able
to recover the image perfectly from an additional $41,808$ complex
noiselet measurements (the total number of real numbers recorded is
$100,000$).

%---------------------------------------------------------------------------------
\section{Proofs}
\label{sec:proof}

\subsection{General strategy}
\label{sec:strategy}

The proof of Theorem~\ref{th:main} follows the program set forth in
\cite{fuchs04sp,candes06ro}.  As detailed in these references, the
signal $x^0$ is the unique solution to \eqref{eq:genl1} if and only if
there exists a {\em dual vector} $\pi\in\R^n$ with the following
properties:
\begin{itemize}
\item $\pi$ is in the row space of $U_\Omega$, 
\item $\pi(t) = \sgn x^0(t)$ for $t\in T$, and
\item $|\pi(t)|<1$ for $t\in T^c$.
\end{itemize}  
We consider the candidate
\begin{equation}
\label{eq:dualvec}
\pi = U_\Omega^* \UOT (\UOT^*\UOT)^{-1}z_0,
\end{equation}
where $z_0$ is a $|T|$-dimensional vector whose entries are the signs
of $x^0$ on $T$, and show that under the conditions in the theorem 1)
$\pi$ is well defined (i.e.\ $\UOT^*\UOT$ is invertible), and given
this 2) $|\pi(t)|<1$ on $T^c$ (we automatically have that $\pi$ is in
the row space of $U_\Omega$ and $\pi(t)=\sgn x(t)$ on $T$).

We want to show that with the support fixed, a dual vector exists with
high probability when selecting $\Omega$ uniformly at random.
Following \cite{candes06ro}, it is enough to show that the desired
properties when $\Omega$ is sampled using a Bernoulli model.  Suppose
$\Omega_1$ of size $m$ is sampled uniformly at random, and $\Omega_2$
is sampled by setting
\[
\Omega_2 := \{k: \delta_k=1\},
\]
where here and below $\delta_1,\delta_2,\ldots,\delta_n$ is a sequence
of independent identically distributed $0/1$ Bernoulli random
variables with
\[
\Prob(\delta_k=1) = m/n.
\]
Then 
\begin{equation}
\label{eq:btou}
\Prob(\mathrm{Failure}(\Omega_1)) \leq 2\Prob(\mathrm{Failure}(\Omega_2))
\end{equation}
(see \cite{candes06ro} for details).  With this established, we will
establish the existence of a dual vector for $x^0$ with high
probability for $\Omega$ sampled using the Bernoulli model.

The matrix $\UOT^*\UOT$ is now a random variable, which can be written as 
\[
\UOT^*\UOT = \sum_{k=1}^n \delta_k u^k \otimes u^k,
\]
where the $u^k$ are the row vectors of $U_T$; $u^k = (U_{t,k})_{t\in
  T}$.

\subsection{Proof of Theorem \ref{teo:largedev}}
\label{sec:largedev}

Our first result, which is an analog to a theorem of Rudelson
\cite[Th.\ 1]{rudelson99ra}, states that if $m$ is large enough, then
on average the matrix $m^{-1}\UtU$ deviates little from the identity.
\begin{theorem}
\label{th:rudelson}
Let $U$ be an orthogonal matrix obeying \eqref{eq:UtU}. Consider a
fixed set $T$ and let $\Omega$ be a random set sampled using the
Bernoulli model. Then
  \begin{equation}
    \label{eq:rudelson}
    \E \| \frac{1}{m} \UOT^*\UOT - I\| 
    \le C_R \cdot \frac{\sqrt{\log |T|}}{\sqrt m} \, 
\max_{1 \le k \le n} \|u^k\|
  \end{equation}
  for some positive constant $C_R$, provided the right-hand side is
  less than 1.  Since the coherence $\mu(U)$ obeys
  \[
\max_{1 \le k \le n}
  \|u^k\| \le \mu(U) \sqrt{|T|}, 
\]
  this implies 
  \begin{equation}
    \label{eq:rudelson2}
     \E \| \frac{1}{m} \UOT^*\UOT - I\| 
    \le C_R \cdot \mu(U) \, \frac{\sqrt{|T| \log |T|}}{\sqrt m}.  
  \end{equation}
\end{theorem}
The probabilistic model is different here than in \cite{rudelson99ra}.
The argument, however, is similar.

\begin{proof}
  We are interested in $\E \|Y\|$ where $Y$ is the random sum
\[
Y = \frac{1}{m} \sum_{k=1}^n \delta_k \,  \ukuk - I. 
\]
Note that since  $U^* U = n I$, 
\[
\E Y = \frac{1}{m} \sum_{k=1}^n \frac{m}{n} \,  \ukuk - I = 
\frac{1}{n} \sum_{k=1}^n \ukuk - I = 0. 
\]

We now use a symmetrization technique to bound the expected value of
the norm of $Y$. We let $Y'$ be an independent copy of $Y$, i.e. 
\begin{equation}
\label{eq:Sp}
Y' = \frac{1}{m} \sum_{k=1}^n  \delta'_k \,  \ukuk - I, 
\end{equation}
where $\delta'_1, \ldots, \delta'_n$ are independent copies of
$\delta_1, \ldots, \delta_n$, and write
\[
\E \|Y\| \le \E \|Y - Y'\|,  
\]
which follows from Jensen's inequality and the law of iterated
expectation (also known as Fubini's theorem). Now let $\epsilon_1,
\ldots, \epsilon_n$ be a sequence of Bernoulli variables taking values
$\pm 1$ with probability 1/2 (and independent of the sequences $\delta$
and $\delta'$). We have 
\begin{align*}
  \E \|Y\| & \le \E_{\delta, \delta'} \|\frac{1}{m} 
\sum_{k=1}^n  (\delta_k - \delta'_k) \,  \ukuk \|\\
  & = \E_\epsilon \, \E_{\delta, \delta'} \|\frac{1}{m} \sum_{1 \le k
    \le n} \epsilon_k (\delta_k -
  \delta'_k) \, \ukuk \|\\
  & \le 2 \E_\epsilon \, \E_{\delta} \|\frac{1}{m} \sum_{1 \le k \le
    n} \epsilon_k \delta_k \, \ukuk \|;
\end{align*}
the first equality follows from the symmetry of the random variable
$(\delta_k - \delta'_k) \, \ukuk$ while the last inequality follows
from the triangle inequality. 

We may know apply Rudelson's powerful lemma \cite{rudelson99ra} which
states that
\begin{equation}
  \label{eq:rudelsonlemma}
  \E_\epsilon  \, \| \sum_{k=1}^n
    \epsilon_k \delta_k \, \ukuk \| \le 
C_R/4 \cdot \sqrt{\log |T|} \cdot \max_{k : \delta_k = 1} 
\|u^k\| \cdot \sqrt{\|\sum_{k=1}^n  \delta_k \ukuk\|}    
\end{equation}
for some universal constant $C_R > 0$ (the notation should make
it clear that the left-hand side is only averaged over
$\epsilon$). Taking expectation over $\delta$ then gives
\begin{align}
\nonumber
\E \|Y\| & \le  C_R/2 \cdot \frac{\sqrt{\log |T|}}{m} \cdot \max_{1 \le k
  \le n} \|u^k\| \cdot \E \sqrt{\|\sum_{k=1}^n  \delta_k \ukuk\|}\\
\label{eq:useful}
& \le  C_R/2 \cdot \frac{\sqrt{\log |T|}}{m} \cdot \max_{1 \le k
  \le n} \|u^k\| \cdot \sqrt{\E \|\sum_{k=1}^n  \delta_k \ukuk\|}, 
\end{align}
where the second inequality uses the fact that for a nonnegative
random variable $Z$, $\E \sqrt{Z} \le \sqrt{\E Z}$. Observe now that 
\[
\E \|\sum_{k=1}^n  \delta_k \ukuk\| = \E \|m Y + m I\| \le m
(\E \|Y\| + 1)
\]
and, therefore, \eqref{eq:useful} gives 
\[
\E \|Y\|  \le a \cdot \sqrt{\E \|Y\| + 1}, \qquad a = 
 C_R/2 \cdot \frac{\sqrt{\log |T|}}{\sqrt m} \cdot \max_{1
  \le k \le n} \|u^k\|.
\]
It then follows that if $a \le 1$,
\[
\E \|Y\|  \le 2a,  
\]
which concludes the proof of the theorem. 
\end{proof}

With Theorem~\ref{th:rudelson} established, we have a bound on the
expected value of $\|m^{-1} \UtU - I\|$.  Theorem~\ref{teo:largedev}
shows that $m^{-1}\UtU$ is close to the identity {\em with high
  probability}, turning the statement about expectation into a
corresponding large deviation result.

The proof of Theorem~\ref{teo:largedev}
uses remarkable estimates about the large deviations of
suprema of sums of independent random variables. Let $Y_1, \ldots,
Y_n$ be a sequence of independent random variables taking values in a
Banach space and let $Z$ be the supremum defined as
\begin{equation}
  \label{eq:Z}
  Z = \sup_{f \in \cF} \, \sum_{i = 1}^n f(Y_i),  
\end{equation}
where $\cF$ is a countable family of real-valued functions.  In a
striking paper, Talagrand \cite{talagrand96ne} proved a concentration
inequality about $Z$ which is stated below, see also
\cite{ledoux01co}[Corollary 7.8].
\begin{theorem}
  \label{teo:LT}
  Assume that $|f| \le B$ for every $f$ in $\cF$, and $\E f(Y_i) = 0$
  for every $f$ in $\cF$ and $i = 1, \ldots, n$. Then for all $t \ge 0$,
  \begin{equation}
    \label{eq:Talagrand}
    \P(|Z - \E Z| > t) \le 3 \exp\left(-\frac{t}{K B} \log \left(1 + \frac{B t}
        {\sigma^2 + B \E \bar Z}\right)\right),
  \end{equation}
  where $\sigma^2 = \sup_{f \in \cF} \, \sum_{i = 1}^n \E f^2(Y_i)$,
  $\bar Z = \sup_{f \in \cF} \, |\sum_{i = 1}^n f(Y_i)|$, and $K$ is a
  numerical constant. 
\end{theorem}
We note that very precise values of the numerical constant $K$ are
known and are small, see \cite{massart00ab} and \cite{rio02in,klein05co}. 

\begin{proof}
{\bf of Theorem~\ref{teo:largedev}}.~~
  Set $Y$ to be the matrix $\frac{1}{m} \UtU - I$ and recall that
  $\frac{1}{n} \sum_{k=1}^n \ukuk = I$, which allows to express
  $Y$ as
\[
Y = \sum_{k=1}^n \left(\delta_k - \frac{m}{n}\right)\,
\frac{\ukuk}{m} : = \sum_{k=1}^n Y_k,
\]
where 
\[
Y_k := \left(\delta_k - \frac{m}{n}\right)\, \frac{\ukuk}{m}. 
\]
Note that $\E Y_k = 0$. We are interested in the spectral norm
$\|Y\|$. By definition,
\[
\|Y\| = \sup_{f_1, f_2} \,\,
\<f_1, Y f_2\> = \sup_{f_1, f_2} \,\,
\sum_{k = 1}^n \<f_1, Y_k f_2\>, 
\]
where the supremum is over a countable collection of unit vectors.
For a fixed pair of unit vectors $(f_1, f_2)$, let $f(Y_k)$ denote the
mapping $\<f_1, Y_k f_2\>$.  Since $\E f(Y_k) = 0$, we can apply
Theorem \ref{teo:LT} with $B$ obeying
\[
|f(Y_k)| \le \frac{|\<f_1, \uk\> \, \<\uk, f_2\>|}{m} \le 
\frac{\|\uk\|^2}{m} \le B, \quad \text{for all } k. 
\] 
As such, we can take $B=\max_{1 \le k \le n} \|\uk\|^2/m$. We now
compute
\begin{align*}
\E f^2(Y_k) & = \frac{m}{n}\left(1-\frac{m}{n}\right) \frac{|\<f_1,
  \uk\> \, \<\uk, f_2\>|^2}{m^2} \\ & \le
\frac{m}{n}\left(1-\frac{m}{n}\right) \, \frac{\|\uk\|^2}{m^2} \,
|\<\uk, f_2\>|^2.
\end{align*}
Since $\sum_{1 \le k \le n} |\<\uk, f_2\>|^2 = n$, we proved that
\[
\sum_{1 \le k \le n} \E f^2(Y_k) \le \left(1-\frac{m}{n}\right)\, 
\frac{1}{m} \, \max_{1\le k \le n} \|\uk\|^2 \le B.
\]
In conclusion, with $Z = \|Y\| = \bar Z$, Theorem \ref{teo:LT} shows
that
\begin{equation}
  \label{eq:almost}
  \P(|~\|Y\| - \E \|Y\|~| > t) ~\le~ 3 \exp\left(-\frac{t}{KB} \log\left(1 +
  \frac{t}{1+\E\|Y\|}\right)\right). 
\end{equation}
Take $m$ large enough so that $\E\|Y\| \le 1/4$ in
\eqref{eq:rudelson2}, and pick $t = 1/4$. Since $B \le \mu^2(U)
|T|/m$, \eqref{eq:almost} gives
\[
P(\|Y\| > 1/2) \le 3e^{-\frac{m}{C_T\mu^2(U) |T|}}, 
\]
for $C_T=4K/\log(6/5)$. Taking $C_1=16C_R$ and $C_2=C_T$ finishes the proof.
\end{proof}

\subsection{Proof of Theorem \ref{th:main}}
\label{sec:proofmain}

With Theorem~\ref{teo:largedev} established, we know that with high
probability the eigenvalues of $\UtU$ will be tightly controlled ---
they are all between $m/2$ and $3m/2$.  Under these conditions, the
inverse of $(\UtU)$ not only exists, but we can guarantee that
$\|(\UtU)^{-1}\|\leq 2/m$, a fact which we will use to show
$|\pi(t)|<1$ for $t\in T^c$.

For a particular $t_0\in T^c$, we can rewrite $\pi(t_0)$ as
\[
\pi(t_0) = \<v^0,(\UtU)^{-1}z\> = \<w^0,z\>, 
\]
where $v^0$ is the row vector of $U_\Omega^*\UOT$ with row index
$t_0$, and $w^0 = (\UtU)^{-1}v^0$.  The following three lemmas give
estimates for the sizes of these vectors. From now on and for
simplicity, we drop the dependence on $U$ in $\mu(U)$.

\begin{lemma}
\label{teo:evo}
The second moment of $Z_0 := \|v^0\|$ obeys
\begin{equation}
  \label{eq:evo}
  \E Z_0^2 \le \mu^2 \, m \, |T|. 
\end{equation}
\end{lemma}
\begin{proof}
  Set $\lambda_k^0 = u_{k,t_0}$. The vector $v^0$ is given by
\[
v^0 = \sum_{k=1}^n \delta_k \, \lambda^0_k \, u_{k} = \sum_{k=1}^n
(\delta_k - \E \delta_k) \, \lambda_k^0 \, \uk, 
\]
where the second equality holds due to the orthogonality of the rows
of $U$: $\sum_{1 \le k \le n} \lambda_k^0 \, u_{k,t} = \sum_{1 \le k
  \le n} u_{k,t_0} \, u_{k,t} = 0$.  We thus can view $v^0$ as a sum
of independent random variables:
\begin{equation}
\label{eq:v0}
v^0 = \sum_{k=1}^n Y_k,\qquad Y_k = (\delta_k-m/n) \lambda_k^0 \uk, 
\end{equation}
where we note that $\E Y_k = 0$. It follows that
\[
\E Z_0^2 = \sum_{k} \E \<Y_k, Y_k\> + \sum_{k' \neq k} \E \<Y_k,
Y_{k'}\> =  \sum_{k} \E \<Y_k, Y_k\>.
\]
Now 
\[ 
\E \|Y_k\|^2 = \frac{m}{n}\left(1-\frac{m}{n}\right) |\lambda_k^0|^2
\|\uk\|^2 \le \frac{m}{n}\left(1-\frac{m}{n}\right) \, |\lambda_k^0|^2
\, \mu^2 |T|.
\]
Since $\sum_k |\lambda_k^0|^2 = n$, we proved that
\[
\E Z_0^2 \le \left(1-\frac{m}{n}\right) \, \mu^2 \, m \, |T|.
\]
This establishes the claim. 
\end{proof}

The next result shows that the tail of $Z_0$ exhibits a Gaussian
behavior.   
\newcommand{\osig}{\overline{\sigma}}
\begin{lemma}
\label{teo:newvo}
  Fix $t_0\in T^c$ and let $Z_0 = \|v^0\|$.  Define $\osig$ as
\[
\osig^2 = \mu^2 \, m \cdot \max(1, \mu |T|/\sqrt{m}).
\]
Fix $a > 0$ obeying $a \le (m/\mu^2)^{1/4}$ if $ \mu |T|/\sqrt{m} > 1$
and $a \le (m/\mu^2 |T|)^{1/2}$ otherwise. Then
\begin{equation}
  \label{eq:great}
  P(Z_0 \ge \mu \sqrt{m |T|} + a \osig) \le e^{-\gamma a^2}, 
\end{equation}
for some positive constant $\gamma > 0$. 
\end{lemma}
The proof of this lemma uses the powerful concentration inequality
\eqref{eq:Talagrand}. 

\begin{proof}
By definition, $Z_0$ is given by 
\[
Z_0 = \sup_{\|f\| = 1} \,\,
\<v^0, f\> = \sup_{\|f\| = 1} \,\,
\sum_{k = 1}^n \<Y_k, f\>
\]
(and observe $Z_0 = \bar Z_0$).  For a fixed unit vector $f$, let
$f(Y_k)$ denote the mapping $\<Y_k, f\>$.  Since $\E f(Y_k) = 0$, we
can apply Theorem \ref{teo:LT} with $B$ obeying
\[
|f(Y_k)| \le |\lambda^0_k| \, |\<f, \uk\>| \le |\lambda^0_k| \,
\|\uk\| \le \mu^2 \, |T|^{1/2} := B.
\]
Before we do this, we also need bounds on $\sigma^2$ and $\E Z_0$. For
the latter, we simply use
\begin{equation}
  \label{eq:trivZ}
  \E Z_0 \le \sqrt{\E Z_0^2} \le \mu \, \sqrt{m \, |T|}.
\end{equation}
For the former 
\[
\E f^2(Y_k) = \frac{m}{n}\left(1-\frac{m}{n}\right) \, |\lambda_k^0|^2\, 
|\<\uk, f\>|^2 \le \frac{m}{n}\left(1-\frac{m}{n}\right) \, \mu^2 \,
|\<\uk, f\>|^2. 
\]
Since $\sum_{1 \le k \le n} |\<\uk, f\>|^2 = n$, we proved that
\[
\sum_{1 \le k \le n} \E f^2(Y_k) \le m \mu^2 \left(1-\frac{m}{n}\right). 
\]
In conclusion, Theorem  \ref{teo:LT} shows that %for $Z = \|S\|$
\begin{equation}
  \label{eq:almost2}
  \P(|Z_0 - \E Z_0| > t) ~\le~ 3 \exp\left(-\frac{t}{KB} \log\left(1 + 
      \frac{Bt}{\mu^2 m + B \mu \sqrt{m |T|}}\right)\right). 
\end{equation}

Suppose now $\osig^2 = B \mu \sqrt{m |T|} \ge \mu^2 m$, and fix $t = a
\osig$. Then it follows from
\eqref{eq:almost2} that 
\[
 \P(|Z_0 - \E Z_0| > t) ~\le~ 3 e^{-\gamma a^2}, 
\]
provided that $B t \le \osig^2$. The same is true if $\osig^2 = \mu^2
m \ge B \mu \sqrt{m |T|}$ and $B t \le \mu^2 m$. We omit the
details. The lemma follows from \eqref{eq:trivZ}.
\end{proof}

\begin{lemma}
  \label{teo:w}
  Let $w^0=(\UtU)^{-1}v^0$. With the same notations and hypotheses as
  in Lemma \ref{teo:newvo}, we have
\begin{equation}
    \label{eq:w}
    \P\left(\sup_{t_0 \in T^c} \|w^0\| \ge 2\mu  \, 
      \sqrt{|T|/m}  + 2 a \osig/m\right) 
    \le n \, e^{-\gamma \, a^2} + 
    \P(\|\UtU\| \le m/2). 
  \end{equation}
\end{lemma}
\begin{proof}
  Let $A$ and $B$ be the events $\{\|\UtU\| \ge m/2\}$ and
  $\{\sup_{t_0 \in T^c} \|v^0\| \le \mu \sqrt{m \, \|T|} + a \,
  \osig\}$ respectively, and observe that Lemma \ref{teo:newvo} gives
  $\P(B^c) \le n e^{-\gamma \, a^2}$. On the event $A \cap B$
\[
\sup_{t_0 \in T^c} \|w^0\| \le \frac{2}{m} \, (\mu \sqrt{m \, |T|} + a \,
  \osig)
\]
The claim follows. 
\end{proof}

\begin{lemma}
  \label{teo:pi}
  Assume that $z(t)$, $t \in T$ is an i.i.d.~sequence of
  symmetric Bernoulli random variables. For each $\lambda > 0$, we have 
  \begin{equation}
    \label{eq:pia}
    \P\left(\sup_{t \in T^c} |\pi(t)| > 1\right) \le 2n e^{-1/2\lambda^2} + 
    \P\left(\sup_{t_0 \in T^c} \|w^0\| > \lambda\right). 
  \end{equation}
\end{lemma}
\begin{proof}
  The proof is essentially an application of Hoeffding's inequality
  \cite{boucheron04co}.  Conditioned on the $w^0$, this inequality
  states that
\begin{equation}
  \label{eq:hoeffding}
  \P\left( |\<w^0,z\>| > 1 \,\, \vert 
\,\, w^0\right) \le 2 e^{-\frac{1}{2\|w^0\|^2}}. 
\end{equation}
Recall that $\pi(t_0) = \<w^0,z\>$. 
It then follows that 
\[
\P \left(\sup_{t_0 \in T^c} |\pi(t_0)| > 1 \,\, | \,\, \sup_{t_0 \in
    T^c} \|w^0\| \le \lambda\right) \le 2n e^{-\frac{1}{2 \lambda^2}},  
\]
which proves the result.
\end{proof}

The pieces are in place to prove Theorem~\ref{th:main}.  Set $\lambda
= 2\mu \, \sqrt{|T|/m} + 2 a \osig/m$.  Combining Lemmas~\ref{teo:pi}
and \ref{teo:w}, we have for each $a>0$ obeying the hypothesis of
Lemma \ref{teo:newvo}
\[
  \Prob\left(\sup_{t\in T^c}|\pi(t)|>1\right)~\leq~ 2n
  e^{-1/2\lambda^2}  + n e^{-\gamma \,
    a^2} + \Prob\left(\|(\UtU)\|\leq m/2\right).  
\]
For the second term to be less than $\delta$, we choose $a$ such that
\[
a^2 = \gamma^{-1} \, \log(n/\delta), 
\]
and assume this value from now on. The first term is less than
$\delta$ if
\begin{equation}
\label{eq:lambda}
\frac{1}{\lambda^2} \ge 2\log(2n/\delta). 
\end{equation}
Suppose $\mu \, |T| \ge \sqrt{m}$. The condition in Lemma \ref{teo:newvo}
is $a \le (m/\mu^2)^{1/4}$ or equivalently
\[
m \ge \mu^2 \, \gamma^{-2} \, [\log(n/\delta)]^2, 
\] 
where $\gamma$ is a numerical constant. In this case, $a\osig \le \mu
\sqrt{m \, |T|}$ which gives
\begin{equation}
\frac{1}{\lambda^2} \ge \frac{1}{16} \frac{m}{\mu^2 \, |T|}. 
\label{eq:cases}
\end{equation}
Suppose now that $\mu \, |T| \le \sqrt{m}$. Then if $|T| \ge a^2$,
$a\osig \le \mu \sqrt{m \, |T|}$ which gives again
\eqref{eq:cases}. On the other hand if $|T| \le a$, $\lambda \le 4 a
\osig/m$ and
\[
\frac{1}{\lambda^2} \ge \frac{1}{16} \frac{m}{a^2 \, \mu^2}. 
\]
To verify \eqref{eq:lambda}, it suffices to take $m$ obeying 
\[
\frac{m}{16 \, \mu^2} \min\left(\frac{1}{|T|}, \frac{1}{a^{2}}\right)
\ge 2\log(2n/\delta).
\]
This analysis shows that the second term is less than $\delta$
if 
\[
m \ge K_1 \, \mu^2 \, \max(|T|, \log(n/\delta)) \, \log(n/\delta)
\]
for some constant $K_1$.  Finally, by Theorem~\ref{teo:largedev}, the
last term will be bounded by $\delta$ if
\[
m \ge K_2 \, \mu^2 \, |T| \, \log(n/\delta)
\]
for some constant $K_2$.  In conclusion, we proved that there exists a
constant $K_3$ such that the reconstruction is exact with probability
at least $1-\delta$ provided that the number of measurements $m$ obeys
\[
m \ge K_3 \, \mu^2 \, \max(|T|, \log(n/\delta)) \, \log(n/\delta). 
\]
The theorem is proved.

%---------------------------------------------------------------------------------
\section{Discussion}

It is possible that a version of Theorem~\ref{th:main} exists that
holds for all sign sequences on a set $T$ simultaneously, i.e.\ we can
remove the condition that the signs are chosen uniformly at random.
Proving such a theorem with the methods above would require showing
that the random vector $w^0=(\UOT^*\UOT)^{-1}v^0$, where $v^0$ is as
in \eqref{eq:v0}, will not be aligned with the fixed sign sequence
$z$.  We conjecture that this is indeed true, but proving such a
statement seems considerably more involved.

The new large-deviation inequality of Theorem~\ref{teo:largedev} can
also be used to sharpen results presented in \cite{candes06qu} about
using $\ell_1$ minimization to find the sparsest decomposition of a
signal in a union of bases.  Consider a signal $f\in\R^n$ that can be
written as a sparse superposition of the columns of a dictionary
$D=(\Psi_1~\Psi_2)$ where each $\Psi_i$ is an orthonormal basis.
In other words $f=Dx^0$, where $x^0 \in\R^{2n}$ has small
support.  Given such an $f$, we attempt to recover $x^0$ by solving
\begin{equation}
\label{eq:bp}
\min_x~\|x\|_{\ell_1}\quad\text{subject to}\quad Dx = f. 
\end{equation}
Combining Theorem~\ref{teo:largedev} with the methods used in
\cite{candes06qu}, we can establish that if
\[
|\supp x| \leq
\mathrm{Const}\cdot\frac{n}{\mu^2(\Psi_1^*\Psi_2)\cdot\log n},
\]
then the following will be true with high probability (where the
support and signs of $x^0$ are drawn at random):
\begin{enumerate}
\item There is no $x \neq x^0$ with $|\supp x|\leq|\supp
  x^0|$ with $f=D x$.  That is, $x^0$ is the sparsest possible
  decomposition of $f$.
\item We can recover $x^0$ from $f$ by solving \eqref{eq:bp}.
\end{enumerate}
This is a significant improvement over the bounds presented in
\cite{candes06qu}, which have logarithmic factors of $(\log n)^6$.

%---------------------------------------------------------------------------------
\small
\bibliographystyle{plain}
\bibliography{PartialMeasurements}

%---------------------------------------------------------------------------------
\end{document}